\documentclass[11pt]{amsart}

\usepackage{mathpazo}
\usepackage{lineno,hyperref}
\modulolinenumbers[5]

\usepackage{amsmath, amssymb}
\usepackage[all]{xy}
\usepackage{calrsfs}
\usepackage{tikz}
\usetikzlibrary{decorations.markings}  %Para usar decorations en tikz.
\usetikzlibrary{snakes}  %Para usar flechas en tikz

\newcommand{\trf}{tr\!f_{\overset{}\xi}}
\newcommand{\Q}{\mathbb Q}
\newcommand{\R}{\mathbb R}
\newcommand{\h}{\mathbb H}
\newcommand{\M}{\mathcal M}
\newcommand{\T}{\mathcal T}
\newcommand{\X}{\mathfrak X}
\newcommand{\Y}{\mathfrak Y}
\newcommand{\Z}{\mathfrak Z}

\DeclareMathOperator{\dif}{Diff}
\DeclareMathOperator{\mcg}{Mod}

\newcommand{\n}{\noindent}

\newtheorem{thm}{Theorem}[section]
\newtheorem{lem}[thm]{Lemma}
\newtheorem{prop}[thm]{Proposition}

\theoremstyle{definition}
\newtheorem{defi}[thm]{Definition}
\newtheorem{rmk}[thm]{Remark}
\newtheorem{ejem}[thm]{Example}

\theoremstyle{proof}

\begin{document}

\nocite{*}

\title{The Nielsen realization problem for \\ non-orientable surfaces}

\author{Nestor Colin}
\address{Departamento de Matem\'aticas, Centro de Investigaci\'on y de Estudios Avanzados \\ del IPN, 
Mexico City 07360, Mexico}
\email{ncolin@math.cinvestav.mx}

\author{Miguel A. Xicot\'encatl} 
\address{Departamento de Matem\'aticas, Centro de Investigaci\'on y de Estudios Avanzados \\ del IPN, 
Mexico City 07360, Mexico}
\email{xico@math.cinvestav.mx}



\maketitle
\markboth{\sc{N. Colin and M. A. Xicot\'encatl}}{\sc{Nielsen Realization Problem for Non-Orientable Surfaces}}

\begin{abstract}
We show the Teichm\"uller space of a non-orientable surface with marked points (considered
as a Klein surface) can be identified with a subspace of the Teichm\"uller space of its orientable double cover. 
Also, it is well known that
 the mapping class group $\mcg(N_g; k)$ of a non-orientable surface can be 
identified with a subgroup of   $\mcg(S_{g-1}; 2k)$,   the mapping class group of its orientable double cover.
These facts together with the classical Nielsen realization theorem are used to prove that every finite subgroup of $\mcg(N_g; k)$ can
be lifted isomorphically to a subgroup of the group of diffeomorphisms $\dif(N_g; k)$. In contrast, we show
the projection  $\dif(N_g)  \to  \mcg(N_g)$ does not admit a section for large $g$.\\
\end{abstract}

%%%%%%%%%%%%%%%%%%%%%%%%%%%%%%%%%%%%%%%%%%%%%
%%%%%%%%%%%%%%%%%%%%%%%%%%%%%%%%%%%%%%%%%%%%%
\section{Introduction}

Let $S_g$ be a closed surface of genus $g$  and $X = \{x_1, \dots, x_k\}$ be a finite subset of $S_g$, of marked points. Let $\dif(S_g ; k )$ denote the group diffeomorphisms of $S_g$ which leave the set $X$ invariant. If $S_g$ is
orientable, let $\dif^{\,+}(S_g ; k)$ be the group of all elements of $\dif(S_g ; k)$ that preserve orientation. The mapping class group 
of $S_g$ with $k$ marked points, is the group
$$
\mcg(S_g ; k) =
\begin{cases}
\dif^{\,+}(S_g; k) / \dif_0(S_g, k)   &  \text{if $S_g$ is orientable}\\
& \\
\dif(S_g; k) / \dif_0(S_g; k) & \text{if $S_g$ is non-orientable},
\end{cases}
$$

\n
where $ \dif_0(S_g; k) $ is the subgroup of $\dif(S_g; k)$ consisting of elements that are isotopic to the identity. If the set of marked points 
is empty, we omit $k$ from the notation and write only $\dif(S_g)$, $\dif^{\,+}(S_g)$, $\dif_0(S_g)$ and $\mcg(S_g)$. 
In 1932, J. Nielsen asked the question about whether finite subgroups of the mapping class group $\mcg(S_g)$ can act on $S_g$
by diffeomorphisms, that is, whether every finite subgroup $G \subset \mcg(S_g)$ can be lifted isomorphically to a subgroup
$\widetilde G \subset \dif(S_g)$.
In the case when $S_g$ is an orientable surface without marked points, this problem was partially solved 
for cyclic groups by Nielsen and other special cases were studied over the years by several authors, see \cite{Zi81}. It is well known that for
$g \geq 2$, the surface $S_g$ admits distinct hyperbolic metrics of constant curvature -1 and in this context,  S. Kerckhoff
 gave a positive answer to the Nielsen realization problem, see \cite{Ke83}.

\begin{thm}[Kerckhoff]
Every finite subgroup of $\mcg(S_g)$ can be realized as a group of isometries of some hyperbolic structure on $S_g$.
\end{thm}

In addition, the mapping class group acts naturally on the Teichm\"uller space $\T(S_g)$ of hyperbolic metrics on $S_g$ and this action is
properly discontinuous. Then Theorem 1.1 is equivalent to

\begin{thm}[Kerckhoff]
Every finite subgroup of $\mcg(S_g)$ acting on $\T(S_g)$ has a fixed point.
\end{thm}

\n
The Teichm\"uller space $\T^k(S_g)$ of genus $g$, $k$ punctured surfaces is defined in a similar way
and the mapping class group
$\mcg(S_g; k)$  acts by Weil-Petersson isometries on $\T^k(S_g)$, see \cite{Wo87}. In this context, 
S. Wolpert gives in \cite{Wo87} a solution of the Nielsen problem with marked points.

\begin{thm}[Wolpert]
Every finite subgroup of $\mcg(S_g; k)$ acting by Weil-Petersson isometries on $\T^k(S_g)$ has a fixed point.
\end{thm}

\n
Moreover, if 
$\mcg^{\, \pm}(S_g; k) = \dif(S_g; k) / \dif_0(S_g; k)$
denotes the extended mapping class group, then Masur and Wolf showed that every group of Weil-Petersson isometries is isomorphic to a subgroup 
of $\mcg^{\,\pm}(S_g; k)$, see \cite{MaWo02}. Thus, Wolpert's result extends the Nielsen realization theorem to the case of finite 
subgroups of $\mcg^{\,\pm}(S_g; k)$.\\

The main purpose of this work is to extend the previous result to the case of non-orientable surfaces $N_g$ with $k$ marked points. Even though this result is known among the experts, to the best of our knowledge there is no proof available 
in the literature. It is well known that the Teichm\"uller space can be defined equivalentely via complex or conformal structures on
Riemann surfaces. In this sense, every non-orientable surface can be endowed with a dianalytic structure or a Klein surface structure and the corresponding Teichm\"uller space is defined as the quotient of all dianalytic structures on $N_g$
modulo the group of diffeomorphisms isotopic to the identity. Using this approach, the Nielsen realization theorem for the
non-orientable case can be deduced from the classical case passing to the orientable double cover. Our main results are the following. \\

\begin{thm}\label{T:main}
Every finite subgroup $G \subset \mcg(N_g; k)$ acting on $\T^k(N_g)$ fixes some point of $\T^k(N_g)$.
\end{thm}

We point out that the Nielsen realization theorem does not hold in general for infinite subgroups of $\mcg(S)$. 
For instance, in the orientable case S. Morita showed that the natural projection   $\dif^{\,+}(S_g)  \to \mcg(S_g)$ does not admit a section for $g\geq 5$, see \cite{Mor87} and \cite{MorBook}. We prove here an analogous statement in non-orientable case.\\

\begin{thm}\label{thm:nonsection}
If $g \geq 35$, the natural projection $p : \dif(N_g)\to \mcg(N_g)$ does not have a section. 
\end{thm}

The paper is organized as follows. In section 2 we review the basics about Klein surfaces and we give a standard model for 
the orientable double cover of $N_g$. In section 3 we recall the relation 
between the mapping 
class groups $\mcg(N_g; k)$ and $\mcg(S_{g-1}; 2k)$
stablished in \cite{HT09} and \cite{GGM18} 
 via the orientable double cover $\pi : S_{g-1} \to N_g$. In section 4
we define the Teichm\"uller space 
$\T^k(N_g)$   using dianalytic structures and show the induced  map 
$\pi^*  :  \T^k(N_g)  \to  \T^{2k}(S_{g-1})$ is injective. In section 5 we use the results from the previous sections to prove Theorem \ref{T:main}.
Finally, in section 6 we use the non-orientable analogues of the Miller-Morita-Mumford classes to prove Theorem 1.5.
Throughout this work we will use the notation $\Sigma = \Sigma^k_g$ when refering to a closed surface of genus $g$, 
with $k$ points removed, where $g$ and $k$ satisfy $2 - g - k < 0$ 
(this condition warranties the surfaces are hyperbolic),
unless indicated otherwise.
To emphasize a surface is orientable we will write $S = S^k_g$ and similarly, we will use $N= N^k_g$
for a non-orientable surface. In some cases we will omit the indices $g$ and $k$, inferring them from the context.
As usual, we shall regard an element $f \in \dif(\Sigma_g; k)$ as an element 
of $\dif(\Sigma^k_g)$ by considering its restriction to the punctured surface $\Sigma^k_g$.
This work is part of the Ph.D. thesis of the first author, written under the supervision of the second author.\\

%%%%%%%%%%%%%%%%%%%%%%%%%%%%%%%%%%%%%%%%%%%%%
%%%%%%%%%%%%%%%%%%%%%%%%%%%%%%%%%%%%%%%%%%%%%
\section{Klein surfaces}

An orientable surface can always be equipped with a complex-analytic structure that makes it a Riemann surface. 
Similarly, a non-orientable surface admits a dianalytic structure giving rise to a Klein surface, a natural 
generalization of a Riemann surface. Since this is not a standard topic, we devote this section to recall the basic 
notions of Klein surfaces, morphisms and the orientable double cover of a Klein surface, which
will be used subsequently to address the Nielsen realization problem in the non-orientable case.
The main references used here are \cite{AG71} and \cite{SS92}. \\

Recall that if $A$ is a nonempty open set in $\mathbb C$ and $f : A \to \mathbb C$ is a map which is differentiable
when regarded as a function of two real variables, then

$$ \tfrac{\partial f}{\partial z}  =  \tfrac{1}{2} \left(\tfrac{\partial f}{\partial x}   -  i \tfrac{\partial f}{\partial y} \right)   
\qquad \text{and} \qquad 
\tfrac{\partial f}{\partial \bar z}  =  \tfrac{1}{2} \left(\tfrac{\partial f}{\partial x}   +  i \tfrac{\partial f}{\partial y} \right).   
$$

We say $f$ is {\em analytic} if  $\tfrac{\partial f}{\partial \bar z} = 0$  on $A$ and {\em antianalytic} 
if  
$\tfrac{\partial f}{\partial z} = 0$ on $A$.
The map $f$ will be called a {\em dianalytic} if its restriction to any connected component of $A$ is analytic or antianalytic.
Let $\Sigma$ be a connected surface without boundary, possibly with a finite number of punctures. 
An atlas $\mathcal U = \{(U_i, \varphi_i)\}_{i\in I}$ of $\Sigma$ will be called a
{\em dianalytic atlas} (resp. {\em analytic atlas}) if all of its transition functions are dianalytic (resp. analytic). 
If the surface has punctures, one requires additionally that every puncture has a neighborhood conformally 
equivalent to the unit disc in $\mathbb C$ minus the origin (without this condition a neighborhood of a puncture could 
be conformally equivalent to a cylinder).
Each pair $(U_i, \varphi_i)$ is called a {\em chart} of $\mathcal U$.
Two dianalytic (resp. analytic) atlases $\mathcal U$, $\mathcal V$ of $\Sigma$ are equivalent if 
$\mathcal U \cup \mathcal V$ is a dianalytic (resp. analytic) atlas as well. An equivalence class $\X$ 
of dianalytic atlases will be called a {\em dianalytic} (resp. {\em analytic}) {\em structure}.
For an oriented smooth surface $\Sigma$ let $\M(\Sigma)$ denote the set of those analytic structures of $\Sigma$ 
which 
{\em agree with the orientation and the differentiable structure}. For a non-orientable surface $\Sigma$ the set 
$\M(\Sigma)$ consists of all dianalytic structures of $\Sigma$ that agree with the differentiable structure.

\begin{defi} A {\em Klein surface} is a topological surface $\Sigma$ together with a dianalytic structure $\X$. 
A {\em Riemann surface} is a topological surface $\Sigma$ together with an analytic structure $\X$.
\end{defi}

\begin{ejem}
Let $\Sigma_g^k$ be a surface of genus $g$ with $k$ punctures. If the surface is hyperbolic, i.e.  
$\chi(\Sigma_g^k)<0$, then its universal cover is the hyperbolic plane  $\h^2$ and $\Sigma_g^k$ is 
the quotient of $\h^2$ by the group $\Gamma$ of covering transformations. Such $\Gamma$ is a non-Euclidean 
crystallographic group and is isomorphic to $\pi_1(\Sigma_g^k)$. Since the group $\Gamma$ does not contain torsion 
elements, the projection induces a dianalytic structure $\X$ on the quotient $\h^2/\Gamma$ and thus $\Sigma_g^k$ 
is a Klein surface. If the $\Sigma_g^k$ is orientable, the group $\Gamma$ is actually a Fuchsian group and one gets 
a Riemann surface.
\end{ejem}

\begin{defi}
A {\em morphism}  of Klein surfaces $ f: (\Sigma, \X)  \to (\Sigma', \Y)$ or a {\em dianalytic map} is a continuous 
map $f : \Sigma \to \Sigma'$ 
such that for all $x \in \Sigma$ there exist dianalytic charts
$(U, \phi)$ and $(V, \psi)$ about $x$ and $f(x)$ 
so that $\psi \circ f \circ \phi^{-1}$ is
a dianalytic mapping on $\phi(U)$. If the surfaces $(\Sigma, \X)$ and $(\Sigma', \Y)$ are Riemann surfaces, the map
$f$ is a {\em morphism} if $\psi \circ f \circ \phi^{-1}$ is analytic when $f$ preserves orientation or antianalytic if $f$
reverses orientation.
\end{defi}

\begin{defi}
For every  $f \in \dif(\Sigma)$ and $\X \in \M(\Sigma)$ one defines the {\em pullback structure} $f^* \X \in \M(\Sigma)$ 
by $f$ as the only structure such that the diffeomorphism $f :  (\Sigma,  f^*\X) \to (\Sigma, \X)$ is a morphism.
\end{defi}

\begin{rmk}
The existence and uniqueness of the pullback $f^*\X \in \M(\Sigma)$ are consequences of \cite{AG71}, 
Theorem 1.5.2.
Notice that when the surface $\Sigma$ is non-orientable, the map $f : (\Sigma, f^*X) \to (\Sigma, \X)$ is dianalytic
and when $\Sigma$ is orientable, 
$f :  (\Sigma, f^*\X) \to (\Sigma, \X)$  is analytic if $f$ is orientation preserving or antianalytic if $f$ is orientation reversing.
Also, the uniqueness of the pullback implies that if $\X, \Y \in \M(\Sigma)$
are two structures such that
$f :  (\Sigma, \Y) \to (\Sigma, \X)$ is a morphism, then $f^*\X = \Y$.
\end{rmk}

%%%%%%%%%%%%%%%%%%%%%%%%%%%%%%%%%%%%%%%%%%%%%
%%%%%%%%%%%%%%%%%%%%%%%%%%%%%%%%%%%%%%%%%%%%%

\medskip

\n
{\bf The orientable double cover.}
Any Klein surface $(\Sigma, \X)$ admits three double covers: the the complex double cover, the orientable double cover 
and the Schottky double cover, although the first two coincide in the case when the surface is connected without 
boundary and a finite number of punctures. We deal here with the orientable double cover, which will be used to 
describe the relation between the Teichm\"uller spaces and the mapping class groups of orientable and non-orientable surfaces.

\begin{defi}
An orientable double cover of a non-orientable Klein surface $(\Sigma, \X)$ is a Riemann surface $(S, \X^0)$ 
together with
a dianalytic map $\pi : S \to \Sigma$ which is an unramified double cover of $\Sigma$, and an antianalytic involution 
$\sigma : S \to S$ such that $\pi \circ \sigma = \pi$.
\end{defi}

The orientable double cover $(S, \pi, \sigma)$ of a Klein surface is unique up to isomorphism of Riemann surfaces: 
if $(S', \pi', \sigma')$ is another orientable double
cover of $\Sigma$, then there is a unique analytic isomorphism $f : S' \to S$ such that $\pi' = \pi \circ f$. 
The existence of the orientable double cover 
%of a {\bf (compact ??)} Klein surface 
can be found in \cite{AG71} Theorem 1.6.7, and the same argument can be carried out in the case of a surface with marked points. 
If $(S, \pi, \sigma)$ is the orientable double cover of a non-orientable Klein surface $N$, then for every
$\X \in \M(N)$ the pullback structure $\pi^*\X \in \M(S)$ is the only structure such that
$\pi : (S, \pi^*\X) \to (N, \X)$ is a dianalytic map. 
Also, it is not hard to see that 
for every $\X \in \M(N)$ the involution $\sigma : (S, \pi^*\X) \to (S, \pi^*\X)$ is antianalytic with respect to $\pi^*\X$.

\begin{rmk}
A pair $(S, \sigma)$ given by a Riemann surface $S$ and a antiholomorphic involution 
$\sigma : S \to S$ is known in the literature as {\em a symmetric Riemann surface}. It is not hard to see that the quotient
$\Sigma = S/\langle \sigma \rangle$ of $S$ under the action of $\sigma$ is a Klein surface. Conversely, 
every Klein surface $\Sigma$ is obtained from a symmetric Riemann surface via the 
orientable double cover $(S, \sigma, \pi)$. Thus, the category of Klein surfaces is equivalent to the category of symmetric Riemann surfaces, see \cite{Braun12} for more details.
\end{rmk}

\bigskip

\n
{\bf Construction of the orientable double cover.}
{\em Case without marked points}.
For $N_g$ a compact connected non-orientable surface of genus $g$, without boundary, the orientable double
cover can be constructed as follows. Let $S_{g-1}$ a closed orientable surface of genus $g-1$, embedded in
$\R^3$ such that $S_{g-1}$ is invariant under reflections in the $xy$-, $yz$-, and $xz$-planes. Let $\sigma : S_{g-1} \to S_{g-1}$ be the orientation
reversing homeomorphism given by
$$ \sigma(x, y, z) =  (-x, -y, -z). $$

Then the quotient  $S_{g-1}/\langle \sigma \rangle$ is homeomorphic to 
$N_g$ and the natural projection $ \pi :  S_{g-1}  \to N_g $
is a (topological) double cover of $N_g$. 
Moreover, if $(N_g, \X)$ is a Klein surface we endow $S_{g-1}$
with the pullback structure $\X^0 = \pi^*\X$. With these structures, the map 
$\pi : S_{g-1} \to N_g$ is dianalytic and the involution $\sigma : S \to S$ is antianalytic, by the above remarks. 
Thus
$(S_{g-1}, \pi, \sigma)$
is an orientable double cover  of the Klein surface $N_g$. \\

\paragraph{Case with marked points}
Let $X = \{x_1, x_2, \dots, x_k\}$ be the set of marked points of $N_g$, and let $\widetilde X = \pi^{-1}(X)$. If $\widetilde x_i \in S_{g-1}$ is such that $\pi(\widetilde x_i) = x_i$ for all $i = 1, \dots, k$,
then

$$  \widetilde X = \{ \widetilde x_1, \sigma(\widetilde x_1),
\widetilde x_2, \sigma(\widetilde x_2), 
\dots, 
\widetilde x_k, \sigma(\widetilde x_k)   \} .$$

Remove the sets of marked points $X$ and $\widetilde X$ from $N_g$ and $S_{g-1}$, respectively, to
obtain punctured surfaces $N_g^k$ and $S_{g-1}^{2k}$. By example 2.2, we can endow the surface $N_g^k$
with a dianalytic structure $\X$. An argument similar to the one in the compact case shows that 
$(S_{g-1}^k, \pi, \sigma)$ is an orientable double cover of the Klein surface $N_g^k$.

\medskip

\begin{figure}[h]\label{Figura1}
\centering

\begin{tikzpicture}[thick, scale=0.8] 
	\tikzstyle{every node}=[font=\footnotesize]

%%%%%%%%%%%%%% Superficie Impar %%%%%%%%%%%%%%%%%%%%%	

	\draw[shading=ball, ball color= cyan] (0,0) ellipse (3.5 cm and 1.5cm);
	\draw (-1.8,0) node {$\ldots$} (1.8,0) node {$\ldots$};
%%Genus
	\foreach \x in {-2.7,-1, 0, 1, 2.7 } {
	\draw [xshift=\x cm,scale=0.6] (-0.65, 0.4) to[out=-80, in=125 ] (-0.5, 0)  (0.5, 0) to[out=75, in=-95 ] (0.55,0.25);
			\draw [xshift=\x cm, fill=white, scale=0.6] (-0.5, 0) .. controls +(0.25, 0.3) and +(-0.25,0.3) .. (0.5, 0) .. controls +(-0.1,-0.4) and +(0.25, -0.4) .. (-0.5, 0);}
%%Genus
%%Brace
	\draw [snake=brace] (-3, 0.4) -- (-0.7,0.4) node[pos=0.5, anchor=south] {$m$};
	\draw [snake=brace, scale=-1] (-3, 0.4) -- (-0.7, 0.4) node[pos=0.5, anchor=north] {$m$};
%%%%%%%%%%%%%%%%%%%%%%%%%%%%%%%%%%%%%%%%%%%

%%%%%%%%%%Flechas de mapeos y texto %%%%%%%%%%%%%%%%
	\draw[->] (0, -2) -- (0, -2.5) node[pos=0.5, anchor=west] {$\pi$};
	\draw[->] (8, -2) -- (8, -2.5) node[pos=0.5, anchor=west] {$\pi$};
	
	\draw (1,-4.2) node[anchor=west] {$S_m \# K \cong$};
	\draw (0.75,-4.8) node[anchor=west] {$N_{2m+2}=N_g$};
	\draw (0.4,-5.4) node[anchor=west] {$K=$ Klein Bottle};
	\draw (9,-4.2) node[anchor=west] { $S_m \# \mathbb{R} P^2 \cong$};
	\draw (8.75,-4.8) node[anchor=west] { $N_{2m+1}=N_g$};	
	\draw (-3.5,2.25) node[anchor=west] { Case:  $  g=2m+2$  };
	\draw (-3.5,1.5) node[anchor=west] { $ S_{g-1}$  };
	\draw (4.5,2.25) node[anchor=west] { Case:  $  g=2m+1$  };
	\draw (4.5,1.5) node[anchor=west] {$ S_{g-1}$  };
%%Punto p y sigma p
	\filldraw [color=red] (-1,1) circle (2pt) node[anchor=north west, color=black] {$p$}; 
	\filldraw [color=red, scale=-1, opacity=0.35] (-1,1) circle (2pt)node[anchor=south east, color=black] { $\sigma(p)$}; 
%%%%%%%%%%%%%%%%%%%%%%%%%%%%%%%%%%%%%%%%%%%

%%%%%%%%%%%%%%%% Superficie Par %%%%%%%%%%%%%%%%%%%
	\begin{scope}[shift={(8,0)}] 
		\draw[shading=ball, ball color= cyan] (0,0) ellipse (3.5 cm and 1.5cm);
		\draw (-1.8,0) node {$\ldots$} (1.8,0) node {$\ldots$};		
		\foreach \x in {-2.7,-1, 1, 2.7 } {
			\draw [xshift=\x cm,scale=0.6] (-0.65, 0.4) to[out=-80, in=125 ] (-0.5, 0)  (0.5, 0) to[out=75, in=-95 ] (0.55,0.25);
			\draw [xshift=\x cm, fill=white, scale=0.6] (-0.5, 0) .. controls +(0.25, 0.3) and +(-0.25,0.3) .. (0.5, 0) .. controls +(-0.1,-0.4) and +(0.25, -0.4) .. (-0.5, 0);}
	%%Braces
	\draw [snake=brace] (-3, 0.4) -- (-0.7,0.4) node[pos=0.5, anchor=south] { $m$};
	\draw [snake=brace, scale=-1] (-3, 0.4) -- (-0.7, 0.4) node[pos=0.5, anchor=north] { $m$};
%%Punto p y sigma p
	\filldraw [color=red] (-1,1) circle (2pt) node[anchor=north west, color=black] {$p$}; 
	\filldraw [color=red, scale=-1, opacity=0.35] (-1,1) circle (2pt)node[anchor=south east, color=black] {$\sigma(p)$}; 
	\end{scope}
%%%%%%%%%%%%%%%%%%%%%%%%%%%%%%%%%%%%%%%%%%%%%

%%%%%%%%%% Cociente de Superficie Par %%%%%%%%%%%%%%
	\begin{scope}[shift={(8,-4.5)}] 
	
		\draw[shading=ball, ball color= cyan] (0,1.5) arc (90:270:3.5 cm and 1.5 cm);
		\draw (-0.3, -1) -- (0.3, 1)  (-0.3, 1) -- (0.3,-1);
		
		\draw[fill=white, very thick] (0,0) ellipse (4 mm  and 1.5 cm) ;
%%Flechass de identificacion%%%%%%%%
		\draw[shift={(-0.4,0)}, fill=black] (0,-3pt) -- (67.5:5pt) -- (112.5:5pt) -- (0,-3pt) -- cycle;
		\draw[shift={(0.4,0)}, fill=black] (0,3pt) -- (67.5:-5pt) -- (112.5:-5pt) -- (0,3pt) -- cycle; 		%%
		
		\draw (-1.8,0) node {$\ldots$}; %genus
		\foreach \x in {-2.7,-1 } {
			\draw [xshift=\x cm,scale=0.6] (-0.65, 0.4) to[out=-80, in=125 ] (-0.5, 0)  (0.5, 0) to[out=75, in=-95 ] (0.55,0.25);
			\draw [xshift=\x cm, fill=white, scale=0.6] (-0.5, 0) .. controls +(0.25, 0.3) and +(-0.25,0.3) .. (0.5, 0) .. controls +(-0.1,-0.4) and +(0.25, -0.4) .. (-0.5, 0);}
%%Braces
		\draw [snake=brace] (-3, 0.4) -- (-0.7,0.4) node[pos=0.5, anchor=south] { $m$};
%%Punto p y sigma p
	\filldraw [color=red] (-1,1) circle (2pt) node[anchor=north, color=black] {\tiny $\pi(p)$}; 
	\end{scope}		
%%%%%%%%%%%%%%%%%%%%%%%%%%%%%%%%%%%%%%%%%%%%%%
	
%%%%%%%%%%%Cociente Superficie Impar %%%%%%%%%%%%%%%%
	\begin{scope}[shift={(0,-4.5)}] 
		\draw [shading=ball, ball color= cyan] (0,1.5) arc (90:270:3.5 cm and 1.5 cm);
		
		\draw (-1.8,0) node {$\ldots$};
		
		\draw[very thick, fill=white] (0,.82) ellipse (3 mm  and 6.7 mm);
		\draw[very thick, fill=white] (0,-.83) ellipse (3 mm  and 6.6 mm);
	%%%Flechas de Identificacion%
		\draw[shift={(-0.3,0.82)}, fill=black] (0,-3pt) -- (67.5:5pt) -- (112.5:5pt) -- (0,-3pt) -- cycle;	
		\draw[shift={(0.3,-0.83)}, fill=black] (0,3pt) -- (67.5:-5pt) -- (112.5:-5pt) -- (0,3pt) -- cycle; 	
		
		\foreach \x in {-2.7,-1 } {
			\draw [xshift=\x cm,scale=0.6] (-0.65, 0.4) to[out=-80, in=125 ] (-0.5, 0)  (0.5, 0) to[out=75, in=-95 ] (0.55,0.25);
			\draw [xshift=\x cm, fill=white, scale=0.6] (-0.5, 0) .. controls +(0.25, 0.3) and +(-0.25,0.3) .. (0.5, 0) .. controls +(-0.1,-0.4) and +(0.25, -0.4) .. (-0.5, 0);}
%%Medio Genero			
			\draw [scale=0.6] (-0.65, 0.4) to[out=-80, in=125 ] (-0.5, 0);			
			\draw [scale=0.6, fill=white]  (0, -0.28) .. controls +(-0,-0) and +(0.25, -0.4) .. (-0.5, 0) .. controls +(0.25, 0.3) and +(0,0) .. (0, 0.25);
%%Braces
		\draw [snake=brace] (-3, 0.4) -- (-0.7,0.4) node[pos=0.5, anchor=south] {$m$};
%%Punto p y sigma p
	\filldraw [color=red] (-1,1) circle (2pt) node[anchor=north, color=black] {\tiny $\pi(p)$}; 
	\end{scope}	
%%%%%%%%%%%%%%%%%%%%%%%%%%%%%%%%%%%%%%%%%%%%%
\end{tikzpicture}

\caption{Orientable double cover of a non-orientable surface $N_g$.}

\end{figure}
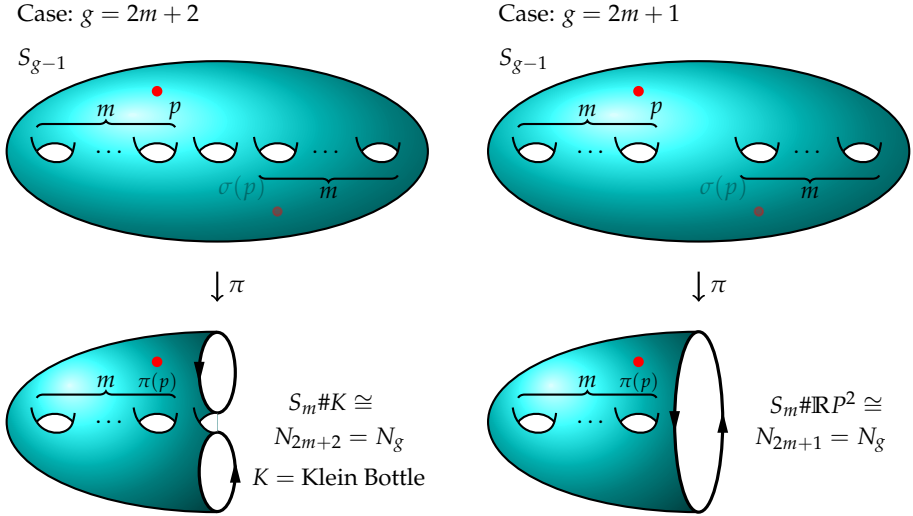

%\medskip
%%%%%%%%%%%%%%%%%%%%%%%%%%%%%%%%%%%%%%%%%%%%%
%%%%%%%%%%%%%%%%%%%%%%%%%%%%%%%%%%%%%%%%%%%%%

\section{Mapping class groups and the orientable double cover}

In this section we review the relationship between the mapping class group of a non-orientable surface $N_g$
and that of its orientable double cover $S_{g-1}$, as studied in \cite{HT09} and \cite{GGM18}. There are several
ways to think of the mapping class group of surface with marked points.
For instance, one may consider a genus $g$, $k$ punctured surface $\Sigma_g^k$ and look at the corresponding
mapping class group $\mcg(\Sigma_g^k)$, where the set of marked points is empty. It is not hard to see that
there is an isomorphism
$$   \mcg(\Sigma_g^k)  \cong  \mcg(\Sigma_g; k) .$$

Thus, we may think of either a surface with marked points or a surface with punctures. The following result is used to
define a natural homomorphism $\phi : \mcg(N_g; k) \to \mcg(S_{g-1}; 2k)$.

\vspace{0.1in}

\begin{prop}
Let $f \in \dif(N_g; k)$ and $\pi : S_{g-1} \to  N_g$ be the orientable double cover. Then $f$ admits exactly two liftings 
$S_{g-1} \to S_{g-1}$, one of which preserves orientation $\widetilde f  \in  \dif^{\,+}(S_{g-1}; 2k)$.
\end{prop}

\n
{\em Proof}. Notice that $\pi' :=  f \circ \pi  :  S_{g-1}  \to  N_g$ is an unramified double cover of $N_g$
and satisfies $\pi'  \circ \sigma = \pi'$. Thus $(S_{g-1}, \pi', \sigma)$ is an orientable double cover
of $N_g$ for suitably chosen structures. By the uniqueness of the orientable double cover, there exists
an analytic isomorphism
$\widehat f  :  S_{g-1} \to S_{g-1}$ such that the following diagram commutes

%\vspace{-0.075in}

$$
\SelectTips{cm}{10}
\xymatrix{
S_{g-1} \ar[r]^{\widehat f} \ar[d]_\pi  \ar[rd]_{\pi'}&  S_{g-1}  \ar[d]^\pi  \\
N_g  \ar[r]_f &  N_g}
$$

\n
Thus the diffeomorphism $\widehat f : S_{g-1} \to S_{g-1}$ is an oriention preserving lifting of $f$ and 
clearly $\sigma\circ \widehat f$ is the only other lift.
Finally, if $X$ is the set of marked points of $N_g$ and $\widetilde X = \pi^{-1}(X)$ is the set of marked points of $S_{g-1}$,
then by the commutativity of the above diagram we have $\widehat f (\widetilde X)=  \widetilde X$ and thus 
$\widehat f \in \dif(S_{g-1} ; 2k)$. \ $\square$ \\

\n
Therefore, there is a natural way to choose a lift of $f \in \dif(N_g; k)$ in a continuous manner by taking
$\widetilde f \in \dif(S_{g-1}; 2k)$ to be orientation preserving. This choice defines a homomorphism
$ \rho  :  \dif(N_g; k)  \to  \dif^{\,+}(S_{g-1}; 2k)$ given by $\rho(f) = \widetilde f$
and induces a homomorphism between the corresponding mapping class groups
$ \phi  :  \mcg(N_g; k)  \to  \mcg(S_{g-1}; 2k) $ such that the following diagram is commutative
$$
\SelectTips{cm}{10}
\xymatrix{
\dif(N_g; k)   \ar[r]^{\rho\quad} \ar[d] & \dif^{\,+}(S_{g-1}; 2k) \ar[d] \\
\mcg(N_g; k)   \ar[r]^{\phi\quad}   &  \mcg(S_{g-1}; 2k)  .\\
}
$$

The following result was proven in \cite{HT09} in the case without marked points, and it was later extended in \cite{GGM18} to the case with marked points. We include it here as it will play an important role in the forthcoming
sections.

\begin{thm}[Hope-Tillmann, Gon\c{c}alves-Guaschi-Maldonado]
Let $N_g$ be a non-orientable surface and $S_{g-1}$ its orientable double cover. Then the following hold:
\begin{enumerate}
\item If $g \geq 3$, the homomorphism $\phi : \mcg(N_g) \to \mcg(S_{g-1})$ is injective.\\

\item If $k \geq 1$, the  homomorphism $\phi : \mcg(N_g; k)  \to  \mcg(S_{g-1}; 2k)$ is injective for all $g \geq 1$.
\end{enumerate}
\end{thm} 
%\hfill $\square$ 

\smallskip
%%%%%%%%%%%%%%%%%%%%%%%%%%%%%%%%%%%%%%%%%%%%%
%%%%%%%%%%%%%%%%%%%%%%%%%%%%%%%%%%%%%%%%%%%%%

\section{The Teichm\"uller space of Klein surfaces}

The Teichm\"uller space can be defined in several equivalent ways, each of which is useful to study specific properties 
of the space. One such definiton is via complex or conformal structures on a Riemann surface. Given that every surface can be equipped 
with a dianalytic structure or a Klein surface structure, it seems natural to extend the definition of the Teichm\"uller space 
as the quotient of the space of all dianalytic structures modulo the group of diffeomorphisms isotopic to the identity. 
The main goal of this section is to establish a relation between the Teichm\"uller spaces of orientable
and non-orientable surfaces by means of the orientable double cover.\\

Let $\Sigma = \Sigma^k_g$ be a genus $g$, $k$ punctured surface without boundary, and let
$\dif_0(\Sigma_g; k)$ denote the subgroup of $\dif(\Sigma_g; k)$ of all diffeomorphisms isotopic to the identity. Two structures 
$\X, \Y \in \M(\Sigma)$ are {\em Teichm\"uller equivalent} if there is a diffeomorphism 
$f \in \dif_0(\Sigma_g; k)$ such that $f : (\Sigma, \X) \to (\Sigma, \Y)$ is a dianalytic map; if $\Sigma$ is a Riemann
surface $f$ is required to be analytic. Equivalently, there exists $f \in \dif_0(\Sigma_g; k)$ such that
$\X =  f^* \Y$.

\begin{defi}
The {\em Teichm\"uller space of a genus $g$, $k$ punctured Klein surface} $\Sigma$, $\T^k(\Sigma_g)$, is the space of 
Teichm\"uller equivalence classes in $\M(\Sigma)$. The elements of the space $\T^k(\Sigma_g)$ will be denoted
by $[\X]$.
\end{defi}

More precisely, the action of $\dif(\Sigma_g; k)$ on $\M(\Sigma)$ by pullbacks restricts to 
an action of $\dif_0(\Sigma_g; k)$ and the space $\T^k(\Sigma_g)$ can be regarded as the quotient
$$   \T^k(\Sigma_g)  =  \M(\Sigma) / \dif_0(\Sigma_g; k) .$$
The set $\M(\Sigma)$ is equipped with a natural topology and we give $\T^k(\Sigma_g)$ the quotient topology.
With this topology $\T^k(\Sigma_g)$ is homeomorphic to $\R^{6g-6+2k}$
in the case of an orientable surface and to $\R^{3g-3+2k}$ in the non-orientable case,
see \cite{PP16}, Theorem 2.2.\\

Consider the orientable double cover $(S, \pi, \sigma)$ of a Klein surface $N$ and recall that for $\X \in \M(N)$ 
the pullback $\pi^*\X$ is the unique structure in $\M(S)$ such that the map $\pi : (S, \pi^*\X) \to (N, \X)$ is dianalytic. Then, the assignment $[\X] \mapsto [\pi^*\X]$ induces a map of Teichm\"uller spaces 
which can be used to identify $\T^k(N_g)$ as a subspace of $\T^{2k}(S_{g-1})$ as follows.

\begin{lem}
Let $(S_{g-1}, \pi, \sigma)$ be the orientable double cover of a non-orientable Klein surface $N_g$. Then the map
induced by the double cover
$$ \pi^* :  \T^k(N_g)  \to  \T^{2k}(S_{g-1}) $$
 is injective.
\end{lem}

\n
{\em Proof.} Let $[\X], [\Y] \in \T^k(N_g)$ be such that $ [\pi^*\X] = [\pi^*\Y]$.
By definition there exists $h \in \dif_0(S_{g-1}; 2k)$ such that    $h : (S_{g-1}, \pi^*\X) \to (S_{g-1}, \pi^*\Y)$
is analytic. Consider the following maps:
$$   (S, \pi^*\X)  \xrightarrow{\;\sigma\;}      
(S, \pi^*\X)  \xrightarrow{\;h\;}   
(S, \pi^*\Y)  \xrightarrow{\;\sigma\;}   
(S, \pi^*\Y)  
$$
Since the involution $\sigma$ is antianalytic with respect to $\pi^*\X$ and $\pi^*\Y$, the 
composition $\sigma \circ h \circ \sigma$ is analytic. On the other hand, 
$  \sigma \circ h \circ \sigma  \simeq  \sigma \circ id \circ \sigma  =  id \simeq h$.
Thus by the uniqueness of the Teichm\"uller Theorems in the
case with marked points (\cite{Hub06}, Corollary 7.2.3), we have $\sigma \circ h \circ \sigma = h$. Therefore
there exists a diffeomorphism $\widehat h : N_g \to N_g$ such that the following diagram commutes
$$
\xymatrix{
(S, \pi^*\X) \ar[r]^h  \ar[d]_\pi  & (S, \pi^*\Y) \ar[d]^\pi   \\
(N, \X)  \ar[r]^{\widehat h}  & (N, \Y)  .
}
$$
Since the maps $\pi\circ h $ and $\pi$ are dianalytic, then
$\widehat h : (N_g, \X) \to (N_g, \Y)$ is also dianalytic (\cite{AG71}, Theorem 1.4.3). It remains to prove that $\widehat h \in \dif_0(N_g; k)$.
Notice the homomorphism $\phi : \mcg(N_g; k) \to \mcg(S_{g-1}; 2k)$, induced by the orientable double cover, satisfies 
$   \phi([\,\widehat h\,]) = [h] =1 .$
By Theorem 3.2 we know that $\phi$ is injective, and so it follows that $\widehat h \in \dif_0(N_g;k)$. 
%We conclude that $[\X] = [\Y]$. 
\hfill$\square$\\

The following result determines the image of $\pi^* :  \T^k(N_g)  \to  \T^{2k}(S_{g-1})$. The proof in the compact
case can be found in \cite{SS92}, Theorem 4.5.1, and we extend it here to the case with marked points.\\

\begin{thm}
Let $(S_{g-1}, \pi, \sigma)$ be the orientable double cover of a non-orien- table Klein surface $N_g$ and 
let $\pi^*:  \T^k(N_g)  \to  \T^{2k}(S_{g-1})$ be
the induced map of Teichm\"uller spaces. Then
$$ \pi^*\left(\T^k(N_g)\right)   =  \left\{ \, [\X]  \in \T^{2k}(S_{g-1})  \; \mid \;  [\sigma^* \X] = [\X] \, \right\} =: \T^{2k}(S_{g-1})^{\sigma^*}$$
\end{thm}

\medskip

\n
{\em Proof.}
If $\Y$ is a dianalytic structure on $N$, then the structure $\pi^*\Y$ on $S$ is such that the involution
$\sigma : (S, \pi^*\Y) \to (S, \pi^*\Y) $ is antianalytic and also
$\sigma^*(\pi^*\Y) = (\pi\circ\sigma)^*\Y= \pi^*\Y$. 
Thus  $\pi^*[\Y] = [\pi^*\Y] \in  \T^{2k}(S_{g-1})^{\sigma^*}$. \\

We proceed to show the reverse inclusion. Let $[\X] \in \T^{2k}(S_{g-1})$ such that $[\sigma^*\X] = [\X]$. 
This means there exists $f \in \dif_0(S)$ such that $f^*\X= \sigma^*\X$ or $\sigma^*(f^*\X) = \X$, 
and thus the following composite is antianalytic
$$  (S, \X)   \xrightarrow{\; \sigma \;} (S, f^*(\X))  \xrightarrow{\; f \;} (S, \X) . $$
Let $\tau := f \circ \sigma : S \to S$ and notice that $\tau^2 \simeq id$. It is actually not hard to see that $\tau^2 =id$. Since  $\tau$ is antianalytic with respect to the structure $\X$, the cyclic group $\langle \tau \rangle$ is a group of automorphisms 
of the Klein surface $(S, \X)$. Thus by \cite{AG71}, Theorem 1.8.4, there exists a unique dianalytic structure $\Y'$ on the surface $N' = S/ \langle \tau \rangle$ such that the canonical 
projection $\pi' : (S, \X) \to (N', \Y')$ is dianalytic.

On the other hand, let us fix a structure $\Z \in \M(N)$ 
and consider the identity map $id  :  (S, \pi^*\Z)  \to  (S, \X)$. By \cite{Hub06}, Corollary 7.2.3, there exists a Teichm\"uller map 
$h :  (S, \pi^*\Z)  \to  (S, \X)$ such that $h \simeq id$ and we consider now the following composite

$$  (S, \pi^*\Z)   \xrightarrow{\; \sigma \;}   
(S, \pi^*\Z)   \xrightarrow{\; h \;}   
(S, \X)   \xrightarrow{\; \tau \;}   
(S, \X)  .
$$
By the remarks in section 2, the involution $\sigma$ is antianalytic with respect to the structure $\pi^* \Z$. If $K_{\tau\circ h \circ \sigma}$ and $K_h$ denote the dilatations of the maps 
$\tau\circ h\circ\sigma$ and $h$, respectively (see \cite{Hub06}), then $K_{\tau\circ h\circ\sigma} = K_h$ since $\tau$
and $\sigma$ are antianalytic. Notice also that
$$ \tau \circ h \circ \sigma \, \simeq  \, \tau \circ \sigma \, = \, f \circ \sigma^2  \, \simeq \, id \, \simeq h ,$$
thus by the uniqueness of the Teichm\"uller Theorems (\cite{Hub06}, Corollary 7.2.3)
we have $\tau \circ h \circ \sigma = h$. Therefore there exists a homeomorphism 
$\widehat h  :  N  \cong  S/\langle \sigma \rangle  \to  N' \cong S/ \langle  \tau \rangle$ such that the following diagram is
commutative:
$$
\xymatrix{
S  \ar[r]^h  \ar[d]_\pi  & S\ar[d]^{\pi'}   \\
N \ar[r]^{\widehat h}  & N'.
}
$$
Define a dianalytic structure by $\Y := (\,\widehat h\,) ^* \,\Y'  \in \M(N)$ and consider the commutative diagram obtained by adding the following structures

$$
\xymatrix{
(S, \pi^*\Y) \ar[r]^h  \ar[d]_\pi  & (S, \X) \ar[d]^{\pi'}   \\
(N, \Y)  \ar[r]^{\widehat h}  & (N', \Y')  .
}
$$
Since the maps $\pi, \pi'$ and $\widehat h$ are dianalytic with respect to the corresponding structures, then by
\cite{AG71}, Theorem 1.4.3, we have that $h  :  (S, \pi^*\Y)  \to  (S, \X)$ is analytic. But $h \simeq id$, thus by definition $\X$  and $\pi^* \Y$ are Teichm\"uller equivalent and therefore
$[\X]  =  [\pi^*\Y] = \pi^*([\Y])$. This completes the proof of the theorem. \hfill $\square$\\

%%%%%%%%%%%%%%%%%%%%%%%%%%%%%%%%%%%%%%%%%%%%%
%%%%%%%%%%%%%%%%%%%%%%%%%%%%%%%%%%%%%%%%%%%%%

\section{Nielsen realization for non-orientable surfaces}

In this section we prove the Nielsen realization theorem in the case of non-orientable surfaces with marked points.
Recall from the previous sections that the orientable double cover $(S_{g-1}, \pi, \sigma)$ of a Klein surface $N_g$, induces injections
$$ \pi^* : \T^k(N_g)   \longrightarrow  \T^{2k}(S_{g-1}) , $$

$$  \phi  :  \mcg(N_g; k)  \to  \mcg(S_{g-1} ;   k)  .$$

\smallskip

\n
We will first show that such maps are compatible with the actions of the mapping class groups on the corresponding
Teichm\"uller spaces. Next we will use the results of Kerckhoff and Wolpert to complete the proof.

\medskip

\begin{lem}
For every $[\X] \in \T^k(N_g)$ and $\alpha \in \mcg(N_g; k)$, the following holds
$$   \pi^*(\alpha \cdot [\X]) =  \phi(\alpha) \cdot\pi^*[\X] .$$
\end{lem}

\medskip

\n
{\em Proof.}
Notice that for  $\alpha = [f]  \in \mcg(N_g; k)$, with  $f \in \dif(N_g; k)$ and  $[\X] \in  \T^k(N_g)$, we have
$$ \pi^*(\alpha \cdot [\X])  =  \pi^*([f] \cdot [\X])  =  [\pi^*(f^*\X)]$$
and
$$  \,\quad \phi(\alpha) \cdot\pi^*[\X]  =  [\rho(f)]\cdot [\pi^*\X]  =  [\rho(f)^*  (\pi^*\X)]  .$$

\smallskip

\n
We will actually show that for every $f \in \dif(N_g;k)$ and $\X \in \M(N_g)$
the following equality holds
$$ \rho(f)^*  (\pi^*\X) =  \pi^*(f^*\X)  .$$
Indeed, by definition of the pullback, the maps
$ \pi : (S, \pi^*\X)  \to  (N,\X)$ and $f : (N, f^*\X)  \to  (N, \X)$ are dianalytic and 
$\rho(f)  :  (S, \rho(f)^*(\pi^*\X))  \to  (S, \pi^*\X)$  is analytic. Consider the following commutative diagram

$$
\xymatrix{
 \!\!(S, \rho(f)^*(\pi^*\X))  \ar[r]^{\quad\rho(f)}    \ar[d]_\pi & (S, \pi^*\X)   \ar[d]^\pi \\
 (N, f^*\X)  \ar[r]^f & (N, \X)
 }
$$

\n
Since the maps $\pi\circ \rho(f)$ and $ f $ are dianalytic, by \cite{AG71}, Theorem 1.4.3, we have   
$  \pi  :  (S,  \rho(f)^*(\pi^*\X))    \to   (N, f^*\X)  $  is dianalytic. Thus by the uniqueness of the pullback we conclude that

$$ \pi^*(f^*\X)   =  \rho(f)^*(\pi^*\X)  ,$$
as we wanted to show.  \hfill $\square$  \\

Next we show the main result.\\

\begin{thm}
Every finite subgroup $G \subset \mcg(N_g; k)$ acting on $\T^k(N_g)$ has a fixed point.
\end{thm}

\n
{\em Proof.}
Recall that the orientable double cover $(S_{g-1}, \pi, \sigma)$ of the surface $N_g$ induces an injective homomorphism

$$ \phi  :  \mcg(N_g; k)  \to  \mcg(S_{g-1}; 2k) , \qquad   [f]  \mapsto [\rho(f)] .$$

\n
Define $H$ to be the subgroup of $\mcg^{\,\pm}(S_{g-1}; 2k)$ generated by $\phi(G)$ and $[\sigma]$ (the isotopy 
class of the involution), and notice this group is isomorphic to $G \times \mathbb Z_2$. By Theorem 1.3
there exists an element $[\Y] \in \T^{2k}(S_{g-1})$ that is fixed under the action of $H$. In particular, for every 
$\alpha \in G$ we have 

\begin{equation}\label{Equation1}
\phi(\alpha)\cdot [\Y] = [\Y]
\end{equation}

\n
and also
$$ [\sigma^* \Y] =[\sigma]\cdot [\Y] = [\Y] . $$

\medskip

\n
By Theorem 4.3 it follows that there exists an element $[\X] \in \T^k(N_g)$ such that 
\begin{equation}\label{Equation2}
\pi^*([\X])  =  [\Y]  .
\end{equation} 

\medskip

\n
We will show that $[\X] \in \T^k(N_g)$ is a fixed point for the action of the group $G$ on the space $\T^k(N_g)$. If 
$\alpha \in G$, then by Lemma 5.1 we have

$$  \pi^*(\alpha \cdot [\X])  =  \phi(\alpha)  \cdot \pi^*([\X])  . $$

\medskip

\n
Thus by equations (\ref{Equation1}) and (\ref{Equation2}) we have

$$ \pi^*(\alpha \cdot [\X])  =  \phi(\alpha)  \cdot  [\Y]  =  [\Y]  =  \pi^*([\X])  .$$

\n
Finally, by Lemma 4.2 we know that $\pi^*$ is injective and therefore

$$ \alpha \cdot  [\X]  =  [\X].   $$
 \hfill $\square$

%%%%%%%%%%%%%%%%%%%%%%%%%%%%%%%%%%%%%%%%%%%%%
%%%%%%%%%%%%%%%%%%%%%%%%%%%%%%%%%%%%%%%%%%%%%

\section{On the non-existence of sections}

Finally, in this last section we prove Theorem 1.5, namely, 
if $g \geq 35$ the projection $p : \dif(N_g)\to \mcg(N_g)$ does not have a section. 
The proof will make use of the non-orientable analogues of the Miller-Morita-Mumford classes, see \cite{Ebert}.
Let $\xi  : E \to B$ be a smooth oriented surface bundle and let $T_v E$ be the vertical tangent bundle. This is a 2-dimensional
real vector bundle on $E$ and thus it has an Euler class   $e \in H^2(E; \mathbb Z)$.
The Miller-Morita-Mumford classes for $\xi$ are defined to be 
$$ \kappa_n(\xi)  :=    \xi_{\overset{}!}\left( e(T_v E)^{n+1}\right)    \in  H^{2n}(B; \mathbb Z) $$

\n
where \, $\xi_{\overset{}!}  :  H^*(E; \mathbb Z)  \to H^{*-2}(B;\mathbb Z)$  \,
 is the {\em umkehr} or the {\em Gysin map}. 
 The generalization of this definition to the non-oriented case uses the Becker-Gottlieb transfer
associated to a smooth surface bundle $\xi : E \to B$. This is a map (in the opposite direction) of suspension spectra 
$\trf :   \Sigma^\infty B_+  \to \Sigma^\infty E_+ $
and thus it induces a morphism in cohomology   $\trf^* :   H^*(E; \mathbb Z)  \to H^*(B; \mathbb Z) $.
The transfer of an oriented surface bundle is related to the umkehr map in the following way: for all $x \in H^*(E; \mathbb Z)$
one has
$$\trf^*(x)  =   \xi_{\overset{}!}   \Big(x\cup  e(T_v E)\Big)  .$$

\n
This implies that
$ \kappa_n(\xi) =  \trf^*\big(e(T_v E)^n\big)$  for the Miller-Morita-Mumford classes of an oriented surface bundle and in particular
$$   \kappa_{2n}(\xi)  =  \trf^* \Big(   p_1(T_v E)^n \Big)      $$
where $p_1(T_v E)$  is the first Pontryagin class of the vertical tangent bundle. This can be generalized to the non-oriented case as follows. For an non-oriented surface bundle  $\eta  :  E \to B$,   N. Wahl in \cite{Wahl} defines classes

$$     \zeta_i(\eta)  :=  tr\!f_{\overset{}\eta}^* \Big(   p_1(T_v E)^i \Big) \,  \in \, H^{4i}(B;\mathbb Z)    $$
where $p_1(T_v E)  \in H^4(E; \mathbb Z)$ is the first Pontryagin class.\\

Recall that for $g \geq 3$, the natural projections  $\dif(N_g)  \to  \mcg(N_g)$ and $\dif^{\,+}(S_{g-1}) \to  \mcg(S_{g-1})$ 
induce homotopy equivalences of classifying spaces by \cite{EE}. On the other hand, the spaces $B \dif(N_g)$ and $B \dif^{\,+}(S_{g-1})$ carry 
universal surface bundles $\eta$ and $\xi$ with fibers $N_g$ and $S_{g-1}$, so the above constructions define cohomology classes   
$$\zeta_i \in  H^{4i}B \dif(N_g) \cong  H^{4i} \mcg(N_g) \text{ \ \ and} $$
$$\kappa_i \in H^{2i} B \dif^{\,+}(S_g) \cong  H^{2i}\mcg(S_{g-1}) $$
which can be used to describe the rational cohomology 
of mapping class groups in the stable range. 
For instance, combining the results of Wahl \cite{Wahl} with those of Galatius, Madsen, Tillmann and Weiss \cite{GMTW}, one gets that
the natural map
$$ \Q[\zeta_1, \zeta_2, \zeta_3, \dots    ]     \longrightarrow     H^*(\mcg(N_g)  ; \Q)$$

\n
is isomorphism in the stable range   $  *   \leq   \frac{g-3}{4}$.
In particular $\zeta_i \neq 0$ in $ H^{4i}( - ; \Q)$  if $g \geq 16 i +3$.
The following result was
proven by Ebert and Randal-Williams in \cite{Ebert}.

\medskip

\begin{thm}[Ebert - Randal-Williams]
Let $\eta : E \to B$ be a non-oriented surface bundle and let $\tilde \eta : \tilde E \to B$ be its fiberwise oriented double cover.
Then the following relation holds for all $n \geq 0$: $$\kappa_{2n}(\tilde\eta) =  2 \cdot \zeta_n(\eta)  .$$
\end{thm}

\medskip

\n
To proceed we will need a couple of lemmas.

\medskip

\begin{lem}
Let $\phi :  \mcg(N_g) \to \mcg(S_{g-1})$ be the morphism induced by the choice of an orientation preserving lifting (see section 3).
 Then for all $i \geq0$:
$$ \phi^*(\kappa_{2i})  =  2 \cdot \zeta_i .$$
\end{lem}

\n
{\em Proof.} Consider the homomorphism  $\rho  :  \dif(N_g) \to  \dif^{\,+}(S_{g-1})$ given by choosing an orientation
preserving lifting and notice that $\dif(N_g)$ acts on $S_{g-1}$ via $\rho$. Then we get a map of bundles
$$
\xymatrix{
S_{g-1}  \ar[d]  \ar@{=}[r] &   S_{g-1}  \ar[d]\\
%\!\!\!\!\!\!\!\!  
E \dif(N_g)   \underset{\dif(N_g)}\times S_{g-1} \ar[r] \ar[d] &  E \dif^{\,+}(S_{g-1}) \underset{\dif^{\,+}(S_{g-1})}\times S_{g-1} \ar[d]  \\
B \dif(N_g)    \ar[r]^{B\rho} &   B \dif^{\,+}(S_{g-1})
}
$$
where the right column is the universal bundle $\xi$ and the left column is the fiberwise orientable double cover $\tilde\eta$
of the universal bundle $\eta$. By naturality of characteristic classes
$$ (B\rho)^*(\kappa_{2i}(\xi)) =  \kappa_{2i}(\tilde\eta)$$
and thus by Theorem 6.2 we have $(B\rho)^*(\kappa_{2i}(\xi))  =  2\cdot \zeta_i(\eta)$. Finally, the commutative 
diagram of groups
$$
\SelectTips{cm}{10}
\xymatrix{
\dif(N_g)   \ar[r]^{\rho\quad} \ar[d] & \dif^{\,+}(S_{g-1}) \ar[d] \\
\mcg(N_g)   \ar[r]^{\phi\quad}   &  \mcg(S_{g-1})  .\\
}
$$
induces a commutative diagram when passing to the cohomology of classifying spaces
$$
\SelectTips{cm}{10}
\xymatrix{
H^*(B \dif(N_g) ; \mathbb Z)  & H^*(B \dif^{\,+}(S_{g-1});  \mathbb Z)   \ar[l]_{B\rho^*}\\
H^*(B \mcg(N_g) ; \mathbb Z)  \ar[u]^\cong &  H^*(B \mcg(S_{g-1}) ; \mathbb Z)  \ar[u]_\cong  \ar[l]_{B\phi^*}   
}
$$
and the result follows.  \hfill $\square$ 

\medskip

\begin{lem}
Let $p : \dif(N_g) \to   \mcg(N_g)$ be the natural projection and let $p^*  : H^*(\mcg(N_g); \Q)  \to H^*(\dif_\delta(N_g) ; \Q)$
be the induced homomorphism, where $\dif_\delta(N_g)$ denotes the group $\dif(N_g)$ considered as a discrete group.
Then for any $i \geq 2$ we have
\; $p^*(\zeta_i) =0  $.
\end{lem}

\n
{\em Proof}.
Notice the commutative diagram
$$
\SelectTips{cm}{10}
\xymatrix{
\dif_\delta(N_g)   \ar[r]^\rho \ar[d]_{p_1} & \dif_\delta^{\,+}(S_{g-1}) \ar[d]^{p_2} \\
\mcg(N_g)   \ar[r]^\phi   &  \mcg(S_{g-1})  .\\
}
$$
gives rise to the following commutative diagram in cohomology

$$
\SelectTips{cm}{10}
\xymatrix{
H^*(\dif_\delta(N_g) ; \Q)  & H^*( \dif_\delta^{\,+}(S_{g-1});  \Q)   \ar[l]_{\rho^*}\\
H^*(\mcg(N_g) ; \Q)  \ar[u]^{p_1^*} &  H^*(\mcg(S_{g-1}) ; \Q)  \ar[u]_{p_2^*}  \ar[l]_{\phi^*}   
}
$$
and thus \ $p_1^*(\phi^*(\kappa_{2i}))  =  \rho^*(p_2^*(\kappa_{2i}))$. By the previous lemma
we have $\phi^*(\kappa_{2i}) = 2 \cdot \zeta_i$ and since $i \geq 2$, Theorem 4.20 of \cite{MorBook} implies that  $p_2^*(\kappa_{2i})=0$.
Thus $p_1^*(2\cdot \zeta_i) = 2 \cdot p_1^*(\zeta_i)=0$.  \hfill $\square$  \\

\medskip

\n
{\em Proof of Theorem \ref{thm:nonsection}}. 
Notice that if $g \geq 35 = 16(2) +3$, the class $\zeta_2 \in H^8(\mcg(N_g); \Q)$
is in the stable range and therefore, is non-zero. Assume that there is a homomorphism  $s : \mcg(N_g) \to  \dif(N_g)$ such that $p\circ s=\rm{id}$
and consider $\dif(N_g)$ as a discrete group.

%\vspace{-0.75in}

$$
\qquad
\xymatrix{
\dif_\delta(N_g)  \ar[d]_p  \\
\mcg(N_g) \ar@/_2pc/[u]_s
}
\qquad \qquad
\xymatrix{
H^*(\dif_\delta(N_g) ; \Q)  \ar@<2ex>[d]^{s^*}  \\
H^*(\mcg(N_g) ;  \Q) \ar@<1ex>[u]^{p^*}
}
$$

\n
Then we have $s^* p^*(\zeta_2) = \zeta_2 \neq 0$. But this contradicts the conclusion $p^*(\zeta_2)=0$ of the previous lemma.
$\square$  

\medskip

\n
\begin{rmk}
In the orientable case when $g=1$, namely for the torus, it is easy to give a section for the projection
$\dif^{\,+}(T^2)  \to   Mod(T^2) = SL(2, \mathbb Z)$ by looking at the natural action of $SL(2, \mathbb Z)$ on $T^2 = \R^2 / \mathbb Z^2$.
In the non-orientable case, when $g=3$, it was shown in \cite{Fico_Mar06} that $Mod(N_3) = GL(2, \mathbb Z)$ and that the standard linear 
action
of $GL(2, \mathbb Z)$ on $T^2$ induces an action on $N_3$, regarded as the blow up of the torus regarded as a real algebraic variety.
This defines
a section for the projection  $ \dif(N_3) \to GL(2, \mathbb Z)$.
Thus, every subgroup of $\mcg(T^2)$ or $\mcg(N_3)$ is realizable as a group of
diffeomorphisms of $T^2$ or $N_3$, regardless if it is finite or not.
\end{rmk}

\medskip

%%%%%%%%%%%%%%%%%%%%%%%%%%%%%%%%%%%%%%%%%%%%%
%%%%%%%%%%%%%%%%%%%%%%%%%%%%%%%%%%%%%%%%%%%%%

\section*{Acknowledgements}
The first author wishes to thank Miguel A. Xicot\'encatl for his support and valuable comments 
while this work was conducted. He also acknowledges the support of CONACYT through the Ph.D. scholarship 
No. 494867.
Both authors are greateful for the financial 
support of CONACYT grant 
CB-2017-2018-A1-S-30345.

\bibliography{References}

\end{document}